\documentclass[12pt, twoside, leqno]{article}

\usepackage{ifthen,latexsym,amssymb,amsmath,amsthm,bbm,fixmath}
\usepackage[shortlabels]{enumitem}
\usepackage{hyperref}

\newtheorem{theorem}{Theorem}[section]
\newtheorem{lemma}[theorem]{Lemma}
\newtheorem{proposition}[theorem]{Proposition}

\newtheorem{corollary}[theorem]{Corollary}

\theoremstyle{definition}

\newtheorem{remark}[theorem]{Remark}

\numberwithin{equation}{section}

\newcommand{\C}[1]{{\protect\mathcal{#1}}}
\newcommand{\B}[1]{{\mathbold #1}}
\newcommand{\I}[1]{{\mathbbm #1}}
\renewcommand{\O}[1]{\overline{#1}}

\newcommand{\ceil}[1]{\lceil #1\rceil}

\newcommand{\floor}[1]{\lfloor #1\rfloor}

\renewcommand{\mid}{:}

\renewcommand{\ldots}{\hspace{0.9pt}.\hspace{0.3pt}.\hspace{0.3pt}.\hspace{1.5pt}}
\renewcommand{\ge}{\geqslant}
\renewcommand{\le}{\leqslant}

\newcommand{\actson}{{\curvearrowright}}
\newcommand{\Matrix}[2]{\left(\begin{array}{#1} #2\end{array}\right)}

\newcommand{\hide}[1]{}

\newcommand{\beq}[1]{\begin{equation}\label{#1}}
	\newcommand{\eeq}{\end{equation}}

\newcommand{\bpf}[1][Proof.]{\smallskip\noindent{\it #1} }
\newcommand{\epf}{\qed \medskip}

\begin{document}
	
	%text specific macros 
	
	\newcommand{\dd}{\,\mathrm{d}}
	\newcommand{\ha}[1]{\cite[#1]{Hassett07itag}}
	\newcommand{\cu}[1]{\cite[#1]{Cutkosky18iag}}
	\newcommand{\gr}[1]{\cite[#1]{Groemer96gafssh}}
    \newcommand{\cgp}[1]{\cite[#1]{ConleyGrebikPikhurko24}}
	\newcommand{\kn}[1]{\cite[#1]{Knapp:lgbi}}
	\newcommand{\TW}[1]{\cite[#1]{TomkowiczWagon:btp}}

	\newcommand{\ga}[2]{(\gamma_{#1})_{#2}}
	\newcommand{\Ad}[1]{\mathrm{Ad}(#1)}
	\newcommand{\Laplace}{\mathop{\mathrm \Delta}}
	\newcommand{\SO}{\mathrm{SO}}
	\newcommand{\QQ}{{\O {\I Q}}}
	\newcommand{\Zero}{{\mathbf 0}}
	\newcommand{\conv}{\mathrm{conv}}
	\newcommand{\Sd}{{\I S^{d-1}}}
	\newcommand{\affdim}{\dim_{\mathrm{aff}}}
	\newcommand{\lcd}{\mathop{\mathrm{lcd}}}
	\newcommand{\Div}{\C D}
	\newcommand{\Borel}{\C B}
	\newcommand{\Frac}{\C F}
	\newcommand{\Id}{1}
	\newcommand{\Ng}{{\C N}}

\newcommand{\ftop}{f^{\mathrm{top}}}
\newcommand{\fmeas}{f^{\mathrm{Haar}}}

	\author{Jan Greb\'\i k\thanks{Supported by MSCA Postdoctoral Fellowships 2022 HORIZON-MSCA-2022-PF-01-01 project BORCA grant agreement number 101105722.}\\
	Computer Science Institute,
Charles University\\
Malostranské nám. 25
118 00 Praha 1, Czechia
	 \and 
		Christian Ikenmeyer\\ Department of Computer Science and Mathematics Institute,\\
		 University of Warwick,
		Coventry CV4~7AL, UK
	\and
		Oleg Pikhurko\thanks{Supported by ERC Advanced Grant 101020255.}\\
		Mathematics Institute and DIMAP,\\
				University of Warwick,
				Coventry CV4 7AL, UK. 
}

\title{Fractional Divisibility of Spheres\\ with Partially Generic Sets of Rotations}

\hide{\classification[%
03E15, %  	Descriptive set theory
33C55, %  	Spherical harmonics
26A21, %  	Classification of real functions; Baire classification of sets and functions
%42C10, %  	Fourier series in special orthogonal functions (Legendre polynomials, Walsh functions, etc.)
43A90, % 	Harmonic analysis and spherical functionsXXxXX
57M60%  	Group actions on manifolds and cell complexes in low dimensions
]{28A05%  	Classes of sets (Borel fields, $\sigma$-rings, etc.), measurable sets, Suslin sets, analytic setsYYyYY
% 54H05%  	Descriptive set theory (topological aspects of Borel, analytic, projective, etc. sets)
}
}

%\keywords{Euclidean sphere, divisibility under a group action, measurable set, special orthogonal group}

\maketitle

\begin{abstract}
 We say that an $r$-tuple $(\gamma_1,\ldots,\gamma_r)$ of special orthogonal $d\times d$ matrices \emph{fractionally divides} the $(d-1)$-dimensional sphere $\I S^{d-1}$ if there is a non-constant function  $f\in L^2(\I S^{d-1})$ such that its translations by $\gamma_1,\ldots,\gamma_r$ sum up to the constant-1 function. Our main result shows, informally speaking, that fractional divisibility is  impossible if at least  $r/2$ rotations are ``generic".  
\end{abstract}

%\begin{lemma}\label{lm:necessary} $f\in L^2(\Sd)$

\section{Introduction}

Let an integer $d\ge 2$ be fixed throughout the paper. Let $\SO(d)$ denote the group of  $d\times d$  special orthogonal real matrices.%
\hide{
\emph{special orthogonal} $d\times d$  matrices, that is, real $d\times d$ matrices $M$ such that the determinant of $M$ is 1 and $M^TM=I_d$, where $I_d$ denotes the identity $d\times d$ matrix.
}
The elements of this group are naturally identified with orientation-preserving 
isometries of the Euclidean unit sphere 
$$
\Sd:=\{\B x\in \I R^d\mid\|\B x\|_2=1\},
$$ 
 %and we will often refer to them as \emph{rotations}. Thus 
 where the action  $\gamma.\B x$ of $\gamma\in \SO(d)$ on the (column) vector $\B x\in\Sd$ is the matrix product $\gamma\B x\in\Sd$.
Also, $\SO(d)$ has the natural left action on the functions on $\Sd$, namely
%Note that the group $\SO(d)$ acts naturally on $L^2(\I S^{d-1},\mu)$ via the \emph{shift action} 
\beq{eq:Action}
(\gamma.f)(\B v):=f(\gamma^{-1}.\B v),\quad \mbox{for }\gamma\in \SO(d),\ f:\I S^{d-1}\to\I R\mbox{ and } \B v\in\I S^{d-1}.
\eeq

Let $\mu$ denote the \emph{spherical measure} on $\I S^{d-1}$, which can be defined as the $(d-1)$-dimensional Hausdorff measure with respect to the geodesic distance
%\emph{arc-length distance} 
on the sphere. This measure is invariant under the action $\SO(d)\actson \Sd$, so the action in~\eqref{eq:Action} preserves the Hilbert space $L^2(\Sd,\mu)$.
% (where the distance between $\B x,\B y\in \I S^{d-1}$ is the angle between the vectors $\B x$ and~$\B y$). 
%We call a subset of $\Sd$ \emph{measurable} if it belongs to the $\mu$-completion of the Borel $\sigma$-algebra on~$\Sd$.

Let $\C D_{d,r}$ be the set of $r$-tuples $\B\gamma=(\gamma_1,\ldots,\gamma_r)\in\SO(d)^r$
such that $\Sd$ is \emph{fractionally $\B\gamma$-divisible} (when we also say that $\B\gamma$ \emph{fractionally divides} $\Sd$), meaning that  there is a function $f\in L^2(\Sd,\mu)$ such that $\sum_{i=1}^r \gamma_i.f$ is equal to $1$ $\mu$-a.e.\ and $f$ differs from $1/r$ on a set of positive measure (or, equivalently, $f$ is not a constant function $\mu$-a.e.). 

This notion was studied in Conley, Greb\'ik and Pikhurko~\cite{ConleyGrebikPikhurko24} who proved in \cite[Theorem~1.1]{ConleyGrebikPikhurko24} that, for any $d,r\ge 2$,  the set $\C D_{d,r}$ cannot contain any \emph{generic} $r$-tuple $\B\gamma\in\SO(d)^r$, where the notion ``generic" is defined in the algebraic geometry sense. Here we do not give the exact definition but point out that the set of non-generic elements of the topological group $\SO(d)^r$ both is meager and has Haar measure 0, see~\cgp{Lemma~1.3}. Thus we get the following corollary.

\begin{corollary}[Conley, Greb\'ik and Pikhurko~\cite{ConleyGrebikPikhurko24}]
\label{cr:cgp} For any integers $d,r\ge 2$, the subset $\C D_{d,r}$ of the compact group $\SO(d)^r$
%of those $r$-tuples in  $\SO(d)^r$ that fractionally divide $\Sd$ 
is both meager and null.\end{corollary}

We say that an $r$-tuple of matrices $(\gamma_1,\ldots,\gamma_r)\in\SO(d)^r$ \emph{measurably divides} $\Sd$ if there is a measurable subset $A$ of $(\Sd,\mu)$ such that its translates $\gamma_1.A,\dots,\gamma_r.A$ (defined in the obvious way) form a partition of $\Sd$, without a single point being omitted or multiply covered. The authors of~\cite{ConleyGrebikPikhurko24} were motivated by the question of Wagon~\cite[Question 4.15]{Wagon:btp} (which is Question 5.15 in the newer edition \cite{TomkowiczWagon:btp} of the book) whether there is a triple of rotations in $\SO(3)$ %$(\gamma_1,\gamma_2,\gamma_3)\in\SO(3)^3$ 
that 
measurably divides $\I S^2$. (The related question whether the Axiom of Choice can be avoided here altogether, not even on a $\mu$-null set, was earlier asked by Mycielski~\cite{Mycielski57,Mycielski57a}.)
The obvious more general version of the question of Wagon, where we consider arbitrary pairs $(r,d)$, has known answer for all $d\not=3$ (see~\cgp{Page 26} for references to these results). However, the cases when $d=3$ and $r\ge 3$ still remain open. Since measurable divisibility is a stronger version of fractional divisibility (namely, we additionally require that $f$ is $\{0,1\}$-valued and $\sum_{i=1}^r \gamma_i.f$ is 1 everywhere), \cite[Theorem~1.1]{ConleyGrebikPikhurko24} implies that any $r$-tuple of matrices achieving it must have entries that satisfy some non-trivial polynomial relation with integer coefficients.

\hide{
Fractional divisibility is of independent interest on its own. For example, Kolountzakis, Lagarias and Wang~\cite{KolountzakisLagarias96,LagariasWang96} initiated systematic study of fractional divisibility of the line $\I R^1$ under the group of translations and there is rich follow-up literature. 
Interestingly, a key role in~\cite{KolountzakisLagarias96,LagariasWang96} was played by periodic sets of translations; this, for a given period, corresponds to fractional division under the action of $\SO(2)$ on the circle~$\I S^1$. 
}

Fractional divisibility of the real line $\mathbb{R}$, better known under the name \emph{(translational) tilings by a function}, has been extensively studied from the perspective of Fourier analysis, see the survey \cite{KolountzakisLev2021}.
A key role in the seminal results~\cite{KolountzakisLagarias96,LagariasWang96} is played by periodic sets of translations; this, for a given period, corresponds to fractional division under the action of $\SO(2)$ on the circle~$\I S^1$.
This special case, when $d=2$, connects our investigation with this area; however, note that when $d>2$ there is a stark difference between commutative and non-commutative divisibility, see for instance the references in \cite[Section~6.4]{KolountzakisLev2021} and the results about translational tilings of tori in \cite{GrebikGreenfeldRozhonTao23imrn}.

The main result of this note (Theorem~\ref{th:main}) is to strengthen Corollary~\ref{cr:cgp} by allowing an adversary to fix some number of matrices. With this in mind, we define $\ftop(d,r)$ (resp.\ $\fmeas(d,r)$) to be the smallest integer $\ell$ such that for every $(r-\ell)$-tuple $\B \gamma'=(\gamma_{\ell+1},\dots,\gamma_{r})\in\SO(d)^{r-\ell}$ the section 
 $$
  \C D_{d,r}(\B \gamma'):=\{(\gamma_1,\dots,\gamma_\ell)\in\SO(d)^\ell\mid  (\gamma_1,\dots,\gamma_r)\in \C D_{d,r}\}
  $$
  (which consists of those $\ell$-tuples $(\gamma_1,\dots,\gamma_\ell)\in\SO(d)^\ell$ such that the extended $r$-tuple $(\gamma_1,\dots,\gamma_r)$ fractionally divides $\Sd$) is a meager (resp.\ null) subset of the group $\SO(d)^\ell$. 
%We expect these two functions to be equal, so let us define $f(d,r):=\max\{\ftop(r,d),\fmeas(d,r)\}$. 
Thus Corollary~\ref{cr:cgp} gives that each of these functions is at most $r$.
Note that, by the Kuratowski-Ulam and Fubini-Tonelli Theorems, if $\ftop(d,r)\le \ell'\le r$ or $\fmeas(d,r)\le \ell'\le r$, then $\ell'$ also satisfies the defining condition.
Thus the following theorem, the main result of this note, strengthens Corollary~\ref{cr:cgp}.

\begin{theorem}\label{th:main} 
	For every $d\ge 2$ and $r\ge 2$,  each of $\ftop(d,r)$ and $\fmeas(d,r)$ is at most $\ell$, where $\ell:=\floor{r/2}$ if $d\ge 3$ and $\ell:=1$ if $d=2$.
	\end{theorem}

Let us provide some lower bounds on these functions.
	
\begin{proposition}\label{pr:lower} It holds $\ftop(d,r),\fmeas(d,r)\ge 1$ for any $d,r\ge 2$. Also, $\ftop(d,4),\fmeas(d,4)\ge 2$ for every odd $d\ge 3$.
\end{proposition}

Thus, by Theorem~\ref{th:main} and Proposition~\ref{pr:lower}, we know the following values:
 \begin{eqnarray*}
 \ftop(2,r)\ =\ \fmeas(2,r)&=&1,\quad \mbox{any $r\ge 2$,}\\
 \ftop(d,3)\ =\ \fmeas(d,3)&=&1, \quad \mbox{any $d\ge 2$,}\\
 \ftop(d,4)\ =\ \fmeas(d,4)&=&2, \quad \mbox{any odd $d\ge 3$.}
 \end{eqnarray*}

There is little empirical data but, if the authors are to guess, it seems plausible that each of $\ftop(d,r)$ and $\fmeas(d,r)$ is at most $2$ for all $d\ge 3$.

\section{Spherical harmonics}\label{se:Harmonics}

Our proofs rely on some basic facts about spherical harmonics. For a detailed introduction to this topic, we refer to the book by Groemer~\cite{Groemer96gafssh}. 
Here, we only state the definitions and results that we need. Recall that an integer $d\ge 2$ is fixed throughout this paper. So the dependence on $d$ is often omitted.

A polynomial $p\in\I R[\B x]$ in $\B x=(x_1,\ldots,x_d)$ is called \emph{harmonic} if $\Laplace p=0$, where $$\Laplace:=\frac{\partial^2}{\partial x_1^2}+\ldots+\frac{\partial^2}{\partial x_d^2}$$ is the \emph{Laplace operator}. A \emph{spherical harmonic} is a function from $\Sd$ to the reals which is the restriction to $\I S^{d-1}$ of a harmonic polynomial on~$\I R^d$. Let $\C H$ be the vector space of all spherical harmonics.
For an integer $n\ge 0$, let $\C H_n\subseteq \C H$ be the linear subspace consisting of all functions $f:\I S^{d-1}\to \I R$ that are the restrictions to $\I S^{d-1}$ of some harmonic polynomial $p$ which is homogeneous of degree~$n$, where we regard the zero polynomial as homogeneous of any degree. By \gr{Lemma~3.1.3}, the homogeneous degree-$n$ polynomial $p$ is uniquely determined by $f\in\C H_n$, so we may switch between these two representations without mention.

In order to state some further properties of spherical harmonics, we need the following definitions. Recall that $\mu$ denotes the spherical measure on~$\I S^{d-1}$. Let the shorthand \emph{a.e.}\ stand for $\mu$-almost everywhere. Let $\sigma_d$ denote  the total measure of the sphere: 
$$
\sigma_d:=\mu(\I S^{d-1})=\frac{2\pi^{d/2}}{\Gamma(d/2)}.
$$
%As $d$ is fixed, the dependence on $d$ is usually not mentioned except for $\sigma_d$ (since $\sigma_{d-1}$ will also appear in some formulas). Also, the shorthand \emph{a.e.}\ stands for $\mu$-almost everywhere.

By \gr{Lemma~1.3.1}, the density $\rho$ of the push-forward of $\mu$ under the projection to any coordinate axis is
\beq{eq:rho}
\rho(t):=\left\{\begin{array}{ll}
	\sigma_{d-1}\,(1-t^2)^{(d-3)/2}, & -1<t<1,\\
	0,& \mbox{otherwise.}
\end{array}\right.
\eeq

An important role is played by the \emph{Gegenbauer polynomials} $(P_0,P_1,\ldots)$ which are obtained from $(1,t,t^2,\ldots)$ by the Gram-Schmidt orthonormalization process on $L^2([-1,1],\rho(t)\dd t)$, except they are normalised to assume value 1 at $t=1$ (instead of being unit vectors in the $L^2$-norm). In the special case $d=3$ (when $\rho$ is the constant function), we get the \emph{Legendre polynomials}.
Of course, the degree of $P_n$ is exactly~$n$.

We are now ready to collect some properties that we will need (with references to the statements in~\cite{Groemer96gafssh} that imply them).

\begin{lemma}\label{lm:P} Let $d\ge 2$. For every integer $n\ge 0$ the following holds.
	\begin{enumerate}[(i),nosep]
		\item \label{it:HnInv} 
		The space $\C H_n$ is invariant under the action of $\SO(d)$ (\gr{Proposition 3.2.4}).
%		\item\label{it:PQ} The polynomial $P_n$ has rational coefficients.
        \item\label{it:Dim} The dimension of $\C H_n$ is 
        $$
        N_n:=\binom{d+n-1}{n}-\binom{d+n-3}{n-2},
        $$ 
        where we agree that $\binom{d+n-3}{n-2}=0$ for $n=0$ or~$1$ (\gr{Theorem 3.1.4}).
		\item\label{it:PHarm} For every $\B v\in\I S^{d-1}$, the function $P_n^{\B v}:\I S^{d-1}\to \I R$, defined by 
		\beq{eq:Pnv}
		P_n^{\B v}(\B x):=P_n(\B v\cdot \B x),\quad \mbox{for }\B x\in\I S^{d-1},
		\eeq
		belongs to $\C H_n$, where $\B v\cdot\B x:=\sum_{i=1}^d v_ix_i$ denotes the scalar product of vectors in  $\I R^d$ (\gr{Theorem 3.3.3}).
		\item\label{it:PBasis} There is a choice of $\B v_1,\ldots,\B v_{N_n}\in\I S^{d-1}$ such that the functions $P_n^{\B v_i}$, $i\in [N_n]$, form a basis of the vector space $\C H_n$ (\gr{Theorem 3.3.14}).
		\item\label{it:PScalar} For every $\B u,\B v\in\I S^{d-1}$, we have $\langle P_n^{\B u},P_n^{\B v}\rangle =\frac{\sigma_d}{N_n} P_n(\B u\cdot \B v)$, where $\langle \cdot,\cdot\rangle$ denotes the scalar product on $L^2(\I S^{d-1},\mu)$ (\gr{Proposition 3.3.6 and Theorem 3.4.1}).
		 \item\label{it:Dec} The subspaces $\C H_0,\C H_1,\ldots$ of $L^2(\I S^{d-1},\mu)$ are pairwise orthogonal
		 (\gr{Theorem 3.2.1}) and the linear span of their union is dense in $L^2(\I S^{d-1},\mu)$ (by \gr{Corollary 3.2.7}).
	\end{enumerate}
\end{lemma}

\section{Proofs}

%We follow notation from~\cite{ConleyGrebikPikhurko24} (included below).

We will need the following characterization of the set $\C D_{d,r}$.

\begin{lemma}\label{lm:char} Let $d,r\ge 2$ and take any $\B\gamma=(\gamma_1,\dots,\gamma_r)\in\SO(d)^r$. For an integer $n\ge 1$, define 
 \beq{eq:LinMap}
  \B L_n(g):=\sum_{s=1}^r \gamma_s.g,\quad \mbox{for $g\in\C H_n$}, 
  %,\quad g\in\C H_n.
  \eeq
  and note that $\B L_n(g)\in\C H_n$ by Lemma~\ref{lm:P}\ref{it:HnInv}.
  Then the $r$-tuple $\B\gamma$ fractionally divides $\Sd$ if and only if there is an integer $n\ge 1$ such that the linear map $\B L_n:\C H_n\to\C H_n$ is not invertible
\end{lemma}

\bpf Note that $\B L_n$ is a linear operator on $\C H_n$. Let $\B L_n^*:\C H_n\to\C H_n$ denote its adjoint operator with respect to the scalar product inherited from $L^2(\Sd,\mu)$. Since the action of $\SO(d)$ preserves the measure $\mu$, we have that $\B L_n^*(g)=\sum_{s=1}^r \gamma_s^{-1}.g$ for $g\in\C H_n$. For a unit vector $\B v\in\I S^{d-1}$, define
\beq{eq:Gv'}
G^{\B v}_{n,\B\gamma}:= \B L_n^*(P_n^{\B v})
%=\sum_{s=1}^r P_n^{\gamma_s^{-1}.\B v}
.
\eeq
 %By Lemma~\ref{lm:P}\ref{it:PHarm}, each function $G^{\B v}_{n,\B\gamma}:\I S^{d-1}\to\I R$, as a linear combination of some spherical harmonics $P_n^{\gamma_s^{-1}.\B v}\in\C H_n$, is itself in $\C H_n$. 
 Note that $G^{\B v}_{n,\B\gamma}=\sum_{s=1}^r P_n^{\gamma_s^{-1}.\B v}$, since it holds
 for all $\gamma\in\SO(d)$ and $\B x\in\Sd$ that
 $$
 \big(\gamma^{-1}.P_n^{\B v}\big)(\B x)=P_n^{\B v}(\gamma.\B x)=P_n(\B v\cdot (\gamma.\B x))=P_n((\gamma^{-1}.\B v)\cdot \B x)=\big(P_n^{\gamma^{-1}.\B v}\big)(\B x)
 %,\quad \mbox{for all $\gamma\in\SO(d)$ and $\B x\in\Sd$}
 .
 $$

Suppose first that there is a non-constant $f\in L^2(\I S^{d-1},\mu)$ such that $\sum_{s=1}^r \gamma_s.f=1$ a.e. By Lemma~\ref{lm:P}\ref{it:Dec},  we have the direct orthogonal sum $\C H =\oplus_{n=0}^\infty \C H_n$. Thus we can uniquely write $f=\sum_{n=0}^\infty F_n$  in $L^2(\I S^{d-1},\mu)$ with $F_n\in\C H_n$ for every~$n\ge 0$. Since the action of $\SO(d)$ preserves each space $\C H_n$ by Lemma~\ref{lm:P}\ref{it:HnInv} as well as the scalar product on $L^2(\I S^{d-1},\mu)$,  we have that
$\gamma.f=\sum_{n=0}^\infty \gamma.F_n$ is the harmonic expansion of $\gamma.f\in L^2(\I S^{d-1},\mu)$.

Since the sum  $\sum_{s=1}^r \gamma_s.f$ is a constant function $1$ a.e.\ and $\C H_0$  is exactly the set of constant functions on $\Sd$, we have
by Lemma~\ref{lm:P}\ref{it:Dec}  and the invariance of the scalar product under $\SO(d)$ that, for every integer $n\ge 1$ and every vector~$\B v\in\I S^{d-1}$,
\begin{eqnarray*}
	0&=&\langle P_n^{\B v},1\rangle\ =\ \langle P_n^{\B v},\gamma_1.f+\ldots+\gamma_r.f\rangle\ =\ \langle P_n^{\B v},\gamma_1.F_n+\ldots+\gamma_r.F_n\rangle\\
	&=& \langle P_n^{\B v},\B L_n(F_n)\rangle\ =\ \langle \B L_n^*(P_n^{\B v}),F_n\rangle \ =\ \langle G^{\B v}_{n,\B\gamma}, F_n\rangle.
%	&=&\langle \gamma_1^{-1}.P_n^{\B v}+\ldots+\gamma_r^{-1}.P_n^{\B v},F_n\rangle\ =\ \langle G^{\B v}_{n,\B\gamma}, F_n\rangle.
\end{eqnarray*}
 %where $G^{\B v}_{n,\B\gamma}$ was defined by~\eqref{eq:Gv}.
% Since the functions $G^{\B v}_{n,\B\gamma}$, $\B v\in\I S^{d-1}$, span the whole space $\C H_n$ by Lemma~\ref{eq:GenericGammas}, we must have that $F_n=0$.

Since $f:\Sd\to\I R$ is not a constant function $\mu$-a.e., there is $n\ge 1$ such that $F_n\not=0$. Thus the functions $G^{\B v}_{n,\B\gamma}\in\C H_n$ for $\B v\in\I S^{d-1}$ do not span the whole space. It follows from Lemma~\ref{lm:P}\ref{it:PBasis} that $\B L_n^*$ is not surjective and thus not invertible; of course, the same applies to~$\B L_n$.

Conversely, suppose that $\B L_n$ is not invertible for some $n\ge 1$. By the finite dimensionality of $\C H_n$ there is a non-zero polynomial $g\in\C H_n$ with $\B L_n(g)=0$. The last identity means that $\sum_{s=1}^r \gamma_s.g=0$. Thus, for example, $f:=\frac1r+g$ shows that $\B\gamma$ fractionally divides $\Sd$.
\epf

\begin{remark}\label{rm:1} Note that the function $f=\frac1r+g$ obtained in the end of the proof of Lemma~\ref{lm:char} is a polynomial so the identity $\sum_{s=1}^r \gamma_s.f=1$ holds everywhere without any exceptions. Also, by taking $f:=\frac1r+cg$ with some constant $c$ satisfying $0<c<\frac1r\,\|g\|_\infty^{-1}$, we can additionally ensure that $0<f<1$ everywhere. Thus, by Lemma~\ref{lm:char}, the definition of fractional divisibility will not be affected if we impose these extra restrictions on $f$.
\end{remark}
\hide{
\begin{theorem}\label{th:main'}
	Let $d\ge 2$ and $r\ge 2$. Let $\gamma_{\ell+1},\dots,\gamma_r\in \SO(d)$ be arbitrary elements, where $\ell:=\ceil{r/2}$ if $d\ge 3$ and $\ell:=1$ if $d=2$. Define $\C D$ to consist of those $\ell$-tuples $(\gamma_1,\dots,\gamma_\ell)\in\SO(d)^\ell$ such that the extended $r$-tuple $(\gamma_1,\dots,\gamma_r)$ fractionally divides the sphere $\Sd$.
	
	Then the set $\C D$, as a subset of the topological group $\SO(d)^\ell$, is both meager and null if $d\ge 3$ and at most countable if $d=2$.
\end{theorem}
}

\bpf[Proof of Theorem~\ref{th:main}.] Recall that we are given arbitrary matrices \[
\gamma_{\ell+1},\dots,\gamma_r\in \SO(d)
\]
and we would like to show that the corresponding section of $\C D_{d,r}$ is ``small''.
Let $Z_n$ consist of  those $\B\gamma=(\gamma_1,\dots,\gamma_\ell)$ in $\SO(d)^\ell$ for which $\B L_n$ is non-invertible, where $\B L_n$ is the linear map from Lemma~\ref{lm:char} defined with respect to the extended $r$-tuple $\B\gamma':=(\gamma_1,\dots,\gamma_r)$.
Since the classes of meager and null subsets are closed under countable unions, 
it is enough by Lemma~\ref{lm:char} to prove that, for every $n$, the set $Z_n$  is both null and meager.

First, suppose that $d\ge 3$. 

By Lemma~\ref{lm:P}\ref{it:PBasis}, choose $\B v_1,\ldots,\B v_{N_n}\in\I S^{d-1}$ such that the functions $P_n^{\B v_i}$, $i\in [N_n]$, form a basis of the vector space $\C H_n$. Note that the linear map $\B L_n$ is non-invertible if and only if
$G^{\B v_i}_{n,\B\gamma'}$ for $i\in [N_n]$ are linearly dependent, where the functions $G^{\B v}_{n,\B\gamma'}:=\B L_n^*(P_n^{\B v})$, for $\B v\in\Sd$, are the same as in the proof of Lemma~\ref{lm:char}.
Consider the $N_n\times N_n$ matrices $L=L(\B\gamma)$ and $M$ with entries
\begin{eqnarray*}
L_{ij}&:=&\frac1{\sigma_{d}} \langle G^{\B v_i}_{n,\B\gamma'},P_n^{\B v_j}\rangle,
%\quad\mbox{for } i,j\in [N_n],
\\
M_{ij}&:=&\frac1{\sigma_{d}} \langle P_n^{\B v_i},P_n^{\B v_j}\rangle,\quad\mbox{for } i,j\in [N_n].
\end{eqnarray*}
Since $M$ is up to a scalar factor the Gram matrix of the basis $\{P_n^{\B v_i}\mid i\in [N_n]\}$ of $\C H_n$, we have that $\det M\not=0$.  Write the vectors $G^{\B v_i}_{n,\B\gamma'}$  in this basis: 
 $$
 (G^{\B v_1}_{n,\B\gamma'},\ldots,G^{\B v_{N_n}}_{n,\B\gamma'})^T=A\,(P_n^{\B v_1},\ldots,P_n^{\B v_{N_n}})^T,
 $$
 for some $N_n\times N_n$ matrix~$A$. Then $L$ is the matrix product $AM$. Thus $\det L$ is zero if and only if $G^{\B v_i}_{n,\B\gamma'}$ for $i\in [N_n]$ are linearly dependent.  By Lemma~\ref{lm:P}\ref{it:PScalar}, we have for every $i,j\in [N_n]$ that
 \begin{eqnarray*}
 L_{ij}&=&\frac1{\sigma_d} \sum_{s=1}^r \langle P_n^{\gamma_s^{-1}.\B v_i},P_n^{\B v_j}\rangle \\
 &=& \frac{1}{N_n}\sum_{s=1}^r P_n((\gamma_s^{-1}.\B v_i)\cdot \B v_j)\ =\ \frac{1}{N_n}\sum_{s=1}^r P_n(\B v_i\cdot (\gamma_s.\B v_j)).
 \end{eqnarray*} 
 For fixed $\B v_1,\ldots,\B v_{N_n}$, this represents each $L_{ij}$ as a polynomial in the $d^2\ell$ entries of the matrices $\gamma_1,\ldots,\gamma_\ell$. Thus $\det L$ is also such a polynomial.
 
Let us see what happens if we let $\gamma_1,\dots,\gamma_{\ell}$ to be equal to $\gamma_{\ell+1}$. The linear map $\B L_n$ 
%from Lemma~\ref{lm:char}. Recall that $\B L_n$ 
sends $g\in\C H_n$ to 
$$
 \sum_{s=1}^r \gamma_s.g=(\ell+1)\gamma_{\ell+1}.g+\sum_{s=\ell+2}^r \gamma_s.g\in\C H_n.
 $$
 Looking at the norm coming from $L^2(\Sd,\mu)\supseteq \C H_n$ and using that the action of $\SO(d)$ on this space preserves the norm, we get by the Triangle Inequality that, for all $g\in \C H_n\setminus\{0\}$,
$$
\| \B L_n(g)\|_2\ge (\ell+1) \| \gamma_{\ell+1}.g\|_2- \sum_{s=\ell+2}^r \| \gamma_s.g\|_2= (2\ell+2-r)\, \|g\|_2>0.
$$
 By the finite dimensionality of $\C H_n$, the linear map $\B L_n$ is invertible. 
Furthermore, note that  the invertible adjoint map $\B L_n^*$ sends the basis vectors $P_n^{\B v_1},\dots,P_n^{\B v_{N_n}}$ to respectively $G^{\B v_1}_{n,\B\gamma'},\dots,G^{\B v_{N_n}}_{n,\B\gamma'}$. Thus the latter vectors are linearly independent (which was observed earlier to be equivalent to $\det L\not=0$). Thus the polynomial $\det L$ is non-zero when evaluated at the constant $\ell$-tuple $(\gamma_{\ell+1},\dots,\gamma_{\ell+1})\in\SO(d)^\ell$. 

We are now ready to show that, for each $n\ge 1$, the set $Z_n$ (that consists of those $(\gamma_1,\dots,\gamma_\ell)\in\SO(d)^\ell$ on which $\det L$ vanishes) is both null and meager. Informally speaking, this is true since the zero set of a non-vanishing polynomial $\det L$ on the irreducible variety $\SO(d)^\ell$ must have strictly smaller dimension than the variety itself. For the sake of completeness, we present a proof of this intuitively clear result (which is obtained by an easy adaptation of the proof of \cgp{Lemma~1.3}).

Clearly, $Z_n$ is a closed subset of of $\SO(d)^\ell$ so, in particular, it is measurable.
First, we prove that  $Z_n$ is a null set with respect to the Haar measure $\nu$ on $\SO(d)^\ell$.
Let us recall how the Haar measure $\nu$ can be constructed for the group $\Gamma:=\SO(d)^\ell$ (and, in fact, for any real Lie group), following the presentation in \kn{Sections VIII.1--2}. Namely, choose some linear basis for the Lie algebra $(\mathfrak{so}(d))^\ell$ viewed as the tangent space $T_{(I_d,\ldots,I_d)}$ at the identity $(I_d,\ldots,I_d)\in\SO(d)^\ell$ and, using the translations of these vectors, turn them into left-invariant vector fields
$X_1,\ldots,X_m$. (Note that the Lie algebra $(\mathfrak{so}(d))^\ell$, that consists of all $\ell$-tuples of skew-symmetric matrices, has dimension $m=\binom{d}{2}\ell$ as a vector space.) For each $\B\gamma\in \Gamma$, let $e_1(\B\gamma),\ldots,e_m(\B\gamma)\in T_{\B \gamma}^*$ be the dual basis to $(X_1(\B\gamma),\ldots,X_m(\B\gamma))$. Then $\omega=e_1\wedge\ldots\wedge e_m$ (the skew-symmetric product) is a smooth $m$ form on $\Gamma$, which is positive and left-invariant and thus
defines a Borel left-invariant non-zero measure on $\Gamma$ (\kn{Theorem 8.21}). By the uniqueness, this has to be a multiple of the Haar measure~$\nu$. In particular, any smooth submanifold of $\Gamma$ of dimension (as a manifold) less than $m$ has zero Haar measure (\kn{Equation~(8.25)}). 

The variety $\SO(d)^\ell$ is irreducible and has dimension $m=\binom{d}{2}\ell$ as a real algebraic variety. This is a well-known fact for $\ell=1$ and the proof of this fact for an arbitrary $\ell$ is carefully spelled out in~\cgp{Lemma 8.2}. 

The set $Z_n\subsetneq \SO(d)^\ell$, as an algebraic variety, has dimension smaller than~$m$ which follows from the definition of the \emph{algebraic dimension} of a variety as the maximum length of a strictly nested chains of non-empty irreducible varieties. (Indeed, a chain for $Z_n$ extends to a larger chain for $\SO(d)^\ell$ by adding $\SO(d)^\ell\supsetneq Z_n$ at the top.) The standard results in the theory of \mbox{(semi-)}\-algebraic sets give 
that every bounded variety of algebraic dimension $k$ in some $\I R^n$ admits a finite partition, with each part being a smooth
semialgebraic manifold in $\I R^n$ of dimension at most $k$, see e.g.\ \cite[Theorem 5.38]{BasuPollackRoy06arag}. 
Apply this result to each of finitely many irreducible components $Z$ of $Z_n$.
The dimension $s$ of each obtained manifold $S$ is at most~$\dim Z$. Indeed, pick a point $\B s\in S$ and the projection from $S$ on some $s$ coordinates which is a homeomorphism around $\B s$.  Observe that these $s$ coordinates are algebraically independent in the function field $\I R(Z)$ (because a non-zero polynomial on $\I R^s$ cannot vanish on a non-empty open set). By \ha{Proposition~7.27},  the algebraic dimension of $Z$ is at least $s$, as claimed. 

Thus we covered $Z_n$ by finitely many manifolds of dimension less than $m$. As it was observed earlier, each such manifold has Haar measure 0 (namely, by~\kn{Equation~(8.25)}). We conclude that the Haar measure of $Z_n$ is indeed zero.

Let us show that $Z_n$ is meager in $\SO(d)^\ell$.  Since $Z_n$ is closed, it is enough to show that the relative interior $U$ of $Z_n\subseteq \SO(d)^\ell$ is empty. Suppose on the contrary that $U\not=\emptyset$. Since the compact group $\SO(d)^\ell$ acts transitively on itself by homeomorphisms, finitely many translates of $U$ cover the whole group. As the Haar measure $\nu$ is invariant under this action, we have that $\nu(U)>0$. However, this contradicts the identity $\nu(Z_n)=0$ that we have already proved.

This proves the theorem for $d\ge 3$.

Finally, suppose that $d=2$. While the above arguments (for $d\ge 3$) also apply for $d=2$, we get in fact the stronger conclusion that each section is countable
by explicitly writing the action of $\SO(2)$ on $L^2(\I S^1,\mu)$. Here, we identify $\SO(2)$ with the circle $\I T:=\I R/2\pi \I Z$ (where we take reals modulo $2\pi$) by corresponding each rotation matrix $\Matrix{ll}{\cos \phi& -\sin\phi\\ \sin \phi& \cos \phi}\in \SO(2)$ to its angle $\phi\in \I T$. Take any integer $n\ge 1$. It is well-known (see e.g.\ \gr{Section 3.1}) that the dimension of $\C H_n$ is $2$ and we can take $c,s\in\C H_n$ as its basis vectors, where $c,s$ are the (unique) homogeneous degree-$n$ polynomials that satisfy $c(\cos \alpha,\sin \alpha):= \cos(n\alpha)$ and $s(\cos \alpha,\sin \alpha):= \sin(n\alpha)$ for every $\alpha\in [0,2\pi)$.
%Then the action of $\phi\in\I T$ on $c$ (resp.\ $s$) gives $\cos(n(x-\phi))=\cos(n\phi)c(x)+\sin(n\phi) s(x)$ (resp.\ $\sin(n(x-\phi))=-\sin(n\phi)c(x)+\cos(n\phi)s(x)$), so its action
Then the action of $\phi\in\I T$ on $\C H_n$ in the basis $(c,s)$ is given by the matrix
 $$
 M_\phi:=\Matrix{cc}{\cos(n\phi)&-\sin(n\phi)\\ \sin(n\phi)& \cos(n\phi)}.
 $$
Indeed, the action of $\phi\in\I T$ on, for example, $c\in\C H_n$ is the function 
 $$\cos(n(x-\phi))=\cos(n\phi)c(x)+\sin(n\phi) s(x),\quad x\in\I T.
 $$
  
 Recall that we are given $\gamma_2,\dots,\gamma_r\in \SO(2)$. 
 Let the linear map $g\mapsto \sum_{s=2}^r \gamma_s.g$ from $\C H_n$ to itself  be given in the basis $(c,s)$ by the matrix $K=(K_{i,j})_{i,j=1}^2$. 
 We vary only $\gamma_1=M_\phi$ for $\phi\in \I T$. Let $x:=\cos(n\phi)$ and $y:=\sin(n\phi)=\pm\sqrt{1-x^2}$. Suppose that $y=\sqrt{1-x^2}$ with the other case being analogous. Consider the linear map $\B L_n$ that sends $g\in \C H_n$ to $\sum_{s=1}^r \gamma_s.
g$.
  As in the case $d\ge 3$, we will be done if we show that the set of those $\phi\in \I T$ for which $\B L_nc$ and $\B L_ns$ are linearly dependent is both null and meager. The matrix of $\B L_n$ in the basis $(c,s)$ is 
  $M_\phi+K$. Its determinant is 
  \begin{eqnarray*}
  \det (M_\phi+K)&=&(x+K_{1,1})(x+K_{2,2})\\
  &-&(-\sqrt{1-x^2}+K_{1,2})(\sqrt{1-x^2}+K_{2,1})\\
  &=&\sqrt{1-x^2}(K_{2,1}-K_{1,2}) + x (K_{1,1}+K_{2,2}) + \det K+1.
  \end{eqnarray*}
  It is enough to show that the set of $x\in [-1,1]$ for which this is 0 is countable, since  the map $\phi\mapsto \cos(n\phi)$ has countable pre-images (in fact, each has at most $2n$ elements). The equation $\det(M_\phi+K)=0$ implies that
  $$
  (K_{1,2}-K_{2,1})^2(1-x^2)= (x(K_{1,1}+K_{2,2})+\det K+1)^2.
  $$
 To satisfy this polynomial identity for for an infinite set of $x$, it must be the case that the coefficients at all powers of $x$ match. By looking at the terms linear in $x$, we conclude that $(K_{1,1}+K_{2,2})(\det K+1)=0$. One of the factors is 0 and either case implies that $K_{1,2}=K_{2,1}$ which in turn implies that $\det K=-1$ and $K_{1,1}+K_{2,2}=0$. Thus $K=\Matrix{cc}{a & b\\ b & -a}$ for some $a,b\in \I R$. Note that $K$ is a linear combination of some matrices $M_\psi$. As each matrix $M_\psi$ has equal diagonal entries and opposite off-diagonal entries, we conclude that $a=b=0$, contradicting $\det K=-1$. Thus the set of $x$ for which $\B L_nc,\B L_ns\in \C H_n$ are linearly dependent is in fact finite and so is its pre-image in $\I T$, as desired. This finishes the proof of the theorem.\epf

%\section{Some positive results}\label{se:examples}
\hide{
Recall that, for integers $d,r\ge 2$, let $\ftop(d,r)$ (resp.\ $\fmeas(d,r)$) be the smallest integer $\ell$ such that for every $\gamma_{\ell+1},\dots,\gamma_{r}\in\SO(d)$ the set $\C D$ that consists of those $\ell$-tuples $(\gamma_1,\dots,\gamma_\ell)\in\SO(d)^\ell$ such that the extended $r$-tuple $(\gamma_1,\dots,\gamma_r)$ fractionally divides the sphere $\Sd$ is a meager (resp.\ null) subset of $\SO(d)^\ell$. Since here we establish some lower bounds which are the same for both functions,
%We expect these two functions to be equal, so 
let us define $f(d,r):=\min\{\ftop(r,d),\fmeas(d,r)\}$.

%Thus T\cgp{Theorem 1.1} gives that each of these functions is well-defined (at most $r$) while Theorem~\ref{th:main} gives that $f(d,r)$ is at most $\ceil{r/2}$ if $d\ge 3$ and at most $1$ if $d=2$.

\begin{lemma}\label{lm:r=4} For every odd $d\ge 3$, it holds that $f(d,4)\ge 2$.\end{lemma}
}

\bpf[Proof of Proposition~\ref{pr:lower}.]
We start by showing that $\ftop(d,r)>0$ and $\fmeas(d,r)>0$ for every $d,r\ge 2$.
By definition, this is the same as showing that we can always find an $r$-tuple of rotations that fractionally divides $\mathbb{S}^{d-1}$.
Given arbitrary $d,r\ge 2$, take $\gamma_1,\dots,\gamma_r\in \SO(d)$ that are the rotations of the $(x_1,x_2)$-plane by angles $\frac{2\pi m}{r}$, for $m\in \{0,\dots,r-1\}$ 
(and fix every other standard basis vector).
Define $f:\Sd\to\{0,1\}$ by $f(\B x):=1$ if there are $\rho>0$ and $\theta\in [0,2\pi/r)$ such that $(x_1,x_2)=(\rho\cos \theta,\rho\sin\theta)$, and let $f(\B x):=0$ otherwise. This $f$ shows that $(\gamma_1,\dots,\gamma_r)\in\C D_{d,r}$, so each of $\ftop(d,r)$ and $\fmeas(d,r)$ is at least~$1$.

Let us consider the case when $d=2m+1\ge 3$ is odd and $r=4$. Let $D(c_1,\dots,c_d)$ denote the diagonal matrix with $c_1,\dots,c_d$ on the diagonal and let $c^{(k)}$ denote the sequence where $c$ is repeated $k$ times. 
\hide{
let us give an example for $d=3$. Let $\gamma_2,\dots,\gamma_3$ be the diagonal matrices $D(1,-1,-1)$, $D(-1,1,-1)$ and $D(-1,-1,1)$. We claim that every $\gamma_1\in\SO(3)$ belongs to $\C D$. Indeed, it is easy to see that $\sum_{s=2}^4 \gamma_s.g=-g$ for every $g\in\C H_1$. Let $\B u\in \I S^2$ be the fixed point of of $\gamma_1$, which exists since $d$ is odd. Then $g$ which sends $\B x$ to $\B u\cdot \B x$ is a non-zero element of $\C H_1$ with $\gamma_1.g=g$ and thus $\sum_{s=1}^4 \gamma_s.g=0$ (and e.g.\ $f=\frac14+\frac14g$ is a non-constant function that shows fractional divisibility).

This example easy extends to show that $f(2m+1,4)\ge 2$: for example, we can let
}

Let $\gamma_2,\gamma_3,\gamma_4$ be the diagonal matrices 
$$
 D((-1)^{(2m-1)},-1,1),\ D((-1)^{(2m-1)},1,-1) \mbox{ and } D((1)^{(2m-1)},-1,-1).
 $$
  We have that their sum is $D((-1)^{(d)})$, the diagonal matrix with all diagonal entries equal to $-1$. We can naturally identify the space $\I R^d$ with $\C H_1$ by corresponding $\B u\in\I R^d$ to the linear map $\B x\mapsto \B x\cdot \B u$. Under this identification the action of  $\gamma\in \SO(d)$ on $\C H_1$ is given by the usual matrix product.
  % by $\gamma^{-1}$. 
  Thus
 $\sum_{s=2}^4 \gamma_s.g=-g$ for every $g\in\C H_1$.

Let us show that the section $\C D_{d,4}(\gamma_2,\gamma_3,\gamma_4)$ is the whole group $\SO(d)$. Take any $\gamma_1\in\SO(d)$. Let $\B u\in \I S^{d-1}$ be a fixed point of $\gamma_1$, which exists since $d$ is odd. Then the function $g:\I S^{d-1}\to\I R$ which sends $\B x$ to $\B u\cdot \B x$ is a non-zero element of $\C H_1$ with $\gamma_1.g=g$ and thus $\sum_{s=1}^4 \gamma_s.g=0$. We can now invoke Lemma~\ref{lm:char} (or just observe that e.g.\ $f=\frac14+\frac14g$ is a non-constant function that shows fractional divisibility).
\epf

\section{Concluding remarks}

Our definition of fractional divisibility allows $f$ to be an arbitrary function in $L^2(\Sd,\mu)$. Of course, by varying this requirements and/or adding constraints such that the values of $f$ are restricted to $[0,1]$, we get potentially different notions. Since the main result of this note (Theorem~\ref{th:main}) is about the non-existence of $f$, we picked the largest natural class of functions for which our proof works. This definition is rather natural; also, it admits an alternative characterization (as given by Lemma~\ref{lm:char}) and is not affected when we add some extra restrictions on $f$ (as discussed in Remark~\ref{rm:1}).

\small

\providecommand{\bysame}{\leavevmode\hbox to3em{\hrulefill}\thinspace}
\providecommand{\MR}{\relax\ifhmode\unskip\space\fi MR }
\providecommand{\MRhref}[2]{%
  \href{http://www.ams.org/mathscinet-getitem?mr=#1}{#2}
}
\providecommand{\href}[2]{#2}

%\small

%\bibliography{bibexport}

%\bibliography{oleg,sets,misc,ramsey,enum,number,posets,sat,ex,matroid,design,random,graph,general,geometry,algorithm,Analysis,limits}

\end{document}